\documentclass[11pt]{article}
\usepackage{amsfonts}
\usepackage{amssymb}
\usepackage{amsmath}
\usepackage{t1enc}
\usepackage[T1]{fontenc}
\usepackage{latexsym}
\usepackage[french,english]{babel}
\newcommand{\cad}{c'est-\`a-dire\ }
\usepackage{theorem}
\newtheorem{enonce}{}[section]
\newtheorem{theo}{Th\'eor\`eme}

\newtheorem{cor}[enonce]{Corollaire}

\newtheorem{prop}[enonce]{Proposition}
\newtheorem{lem}[enonce]{Lemme}
\newtheorem{lemme}[enonce]{Lemme}
\newtheorem{defi}[enonce]{Definition}
\newcommand{\pre}{\noindent {\bf Preuve.\ }}

\newcommand{\ff}{\mathcal F}
\newcommand{\hop}{\vskip .3cm\noindent}
\newcommand{\hip}{\vskip .1cm\noindent}

\title{Feuilletages totalement g\'eod\'esiques, flots riemanniens et vari\'et\'es de Seifert.}
\author{Pierre Mounoud\thanks{Universit\'e d'Avignon, laboratoire d'analyse non lin\'eaire et g\'eom\'etrie, 33 rue Louis Pasteur, 84000 Avignon, France.\quad Courriel : {\tt pierre.mounoud@univ-avignon.fr}}\thanks{{\sf Classification math\'ematique par sujet (2000)} : 57R30 ; 53C50}}
\date{}
\begin{document}
\maketitle
\selectlanguage{french}
\begin{abstract}
Nous \'etudions les feuilletages lisses totalement g\'eod\'esiques de codimension $1$  des vari\'et\'es lorentziennes. Nous nous int\'eressons notamment aux relations entre les flots riemanniens et les feuilletages g\'eod\'esiques. Nous prouvons que, quitte \`a prendre un 2-rev\^etement, tout fibr\'e de Seifert poss\`ede un tel feuilletage.
\end{abstract}
\selectlanguage{english}
\begin{abstract}
We study totally geodesic codimension $1$ smooth foliations on Lorentzian manifold. We are in particular interested by the relations between riemannian flows and geodesic foliations. We prove that, up to a $2$-cover, any Seifert bundle admit such a foliation.
\end{abstract}
\selectlanguage{french}
\tableofcontents
\section{Introduction.}
Dans cet article, nous nous sommes int\'er\'ess\'es aux feuilletages totalement g\'eod\'esiques lisses de codimension $1$ des vari\'et\'es lorentziennes. 
Ces feuilletages ont \'et\'e \'etudi\'es auparavant par A. Zeghib dans \cite{Zeg-feuil} dans le cas o\`u les feuilles sont de type lumi\`ere, par C. Boubel, C. Tarquini et l'auteur dans \cite{BMT} dans le cas o\`u les feuilles sont toutes de type temps et par K. Yokumoto dans \cite{yokumoto}  dans le cas mixte c'est-\`a-dire avec des feuilles de diff\'erents types. Remarquons aussi que lorsque les feuilles sont de type espace, ces feuilletages sont exactement les feuilletages totalement g\'eod\'esiques riemanniens (cf. proposition \ref{debut}) qui, dans le cas des vari\'et\'es compl\`etes de dimension $3$ ont \'et\'e \'etudi\'es par Y. Carri\`ere et E. Ghys dans \cite{Car-Ghys}. Notons aussi que les r\'esultats de ces trois articles concernent presque exclusivement les vari\'et\'es compactes de dimension $3$.

Le but de cet article est de voir quelles sont les vari\'et\'es de Seifert ou  les fibr\'es en cercles qui admettent de tels feuilletages et ceci sous les diverses conditions que l'on imposera. Ainsi dans la partie \ref{ndg}, on rappelle ce qu'il en est, en dimension $3$, si on impose aux feuilles d'\^etre toutes non d\'eg\'en\'er\'ees. Dans la partie \ref{lux}, on s'int\'eresse aux feuilletages totalement g\'eod\'esiques de type lumi\`ere des vari\'et\'es de dimension $3$.
 Nous faisons alors la remarque suivante (bien connue) : un feuilletage $\ff$ de type lumi\`ere et de codimension $1$ contient son orthogonal. On en d\'eduit que la classe d'Euler du fibr\'e tangent \`a $\ff$ est  nulle. Cela nous permet de montrer, avec la classification des vari\'et\'es de dimension $3$ poss\'edant un feuilletage lisse dont les feuilles sont des plans des cylindres et des tores (th\'eor\`eme \ref{gloubi}), le principal r\'esultat de cette section~:
\hop
{\bf Th\'eor\`eme \ref{lumos}} {\sl Si un fibr\'e en cercle orientable $M$ sur une surface orientable $\Sigma$ de genre $g>1$ compacte et de nombre d'Euler $\mathrm{eul}(M)$ poss\`ede un 
feuilletage de codimension $1$ lisse totalement g\'eod\'esique de type lumi\`ere 
alors $|\mathrm{eul}(M)|= 2g-2$.}
\hop
Finalement dans la partie \ref{mix}, on regarde le cas g\'en\'eral avec des feuilles de tout types.
Nous nous restreignons en fait au cas o\`u les feuilles de type lumi\`ere sont compactes et en nombre fini.
Tout d'abord nous cherchons  \`a r\'epondre \`a la question~:\\
 <<\'etant donn\'e un feuilletage $\ff$ de codimension $1$  et un flot riemannien $\phi$ quelle position relative de $\ff$ et de $\phi$ permet de conclure qu'il existe une m\'etrique lorentzienne $g$ telle que $\ff$ est totalement g\'eod\'esique et que $\ff$ et $\phi$ sont orthogonaux pour $g$ ?>> 
\vskip .1cm
Pr\'ecisons  que la pr\'esence d'un feuilletage totalement g\'eod\'esique de codimension $1$ n'implique pas forc\'ement celle d'un flot riemannien sur toute la vari\'et\'e (mais voir proposition \ref{critere}).
\vskip .1cm
Naturellement pour chaque feuille $F$ de $\ff$, $\phi$ devra \^etre soit partout transverse soit partout tangent \`a $F$, selon que la feuille devienne de type lumi\`ere ou non. Cependant cette condition est loin  d'\^etre suffisante.
Le th\'eor\`eme \ref{lapinot} r\'epond \`a la question, avec la restriction indiqu\'ee sur les feuilles de type lumi\`ere, et donne des conditions n\'ecessaires et suffisantes. 
Ces conditions portent sur les positions relatives de $\ff$ et $\phi$ et sont assez techniques. Nous renvoyons le lecteur \`a la section \ref{mix} pour un \'enonc\'e exact de ce r\'esultat. Nous donnons tout de m\^eme ici une id\'ee de celles-ci.  D'une part il y a  une hypoth\`ese d'invariance faible de $\ff$ dans la direction de $\phi$ (cf. d\'efinition \ref{adapte} et proposition \ref{voisinage}) au voisinage des feuilles de tangences entre les deux feuilletages.
D'autre part il y a une hypoth\`ese reliant les propri\'et\'es transverses de $\phi$ et la r\'epartition des feuilles de tangences de $\ff$ (cf. d\'efinition \ref{C} et proposition \ref{globi}). En particulier le flot $\phi$ devra \^etre lorentzien sur certains ouverts de la vari\'et\'e.

Malgr\`es les difficult\'es et les restrictions, le th\'eor\`eme \ref{lapinot} permet de r\'epondre \`a notre question initiale et nous prouvons~:
\hop
{\bf Th\'eor\`eme \ref{seifert}} {\sl Quitte \`a prendre un rev\^etement \`a $2$ feuillets, tout fibr\'e de Seifert de dimension $3$ poss\`ede un feuilletage totalement g\'eod\'esique.}
\hop
Ceci comprend les exemples construit par  K. Yokumoto  dans \cite{yokumoto}.
En fait nous montrons que toute vari\'et\'e de dimension $3$ poss\'edant un flot riemannien transversalement orientale poss\`ede un feuilletage totalement g\'eod\'esique. En dimension sup\'erieure la situation se complique
la d\'emonstration en dimension $3$ utilise de fa\c con essentielle le fait que les fibres singuli\`eres sont isol\'ees ce qui n'est plus vrai en dimension plus grande. On obtient toutefois~:
\hop
{\bf Th\'eor\`eme \ref{cercle}} {\sl Tout fibr\'e en cercle orientable poss\`ede un feuilletage totalement g\'eod\'esique mixte de codimension $1$}\hop
En particulier les sph\`eres de dimension impaire poss\`edent donc des feuilletages totalement g\'eod\'esiques.
\section{D\'efinitions}
Dans cette section nous allons d\'efinir certains termes que nous utiliserons dans la suite et que nous avons class\'es par domaine. Pr\'ecisons que les vari\'et\'es consid\'er\'ees dans la suite seront toujours (sauf pr\'ecision contraire) connexes et sans bord.
\subsection{G\'eom\'etrie lorentzienne.}
\begin{defi}
Soit $(M,g)$ une vari\'et\'e lorentzienne. Un sous-espace vectoriel non nul $E$ du tangent en un point $x$ de $M$, not\'e $T_xM$, sera dit de type espace (resp. temps) (resp. lumi\`ere) si la restriction de $g$ \`a $E$ est riemannienne (resp. lorentzienne ou d\'efinie n\'egative) (resp. d\'eg\'en\'er\'ee).

Par extension, on dira qu'une sous-vari\'et\'e $S$ de $M$ est de type espace (resp. temps ou lumi\`ere) si en tout point $x$ de $S$, on a $T_xS$ est de type espace (resp. temps ou lumi\`ere).
\end{defi}
\noindent D\'efinissons maintenant ce qui constitue le coeur de cet article.
\begin{defi}
Soit $(M,g)$ une vari\'et\'e (pseudo-)riemannienne.
Une sous-vari\'et\'e $S$ de $M$ est dite \emph{totalement g\'eod\'esique} si toute g\'eod\'esique de $g$ tangente \`a $S$ en un point est tangente \`a $S$ partout.

Nous dirons qu'un feuilletage $\ff$ de $M$ est totalement g\'eod\'esique (pour la m\'etrique $g$) si toute feuille de $\ff$ est une sous-vari\'et\'e totalement g\'eod\'esique de $(M,g)$.
\end{defi}
Si on consid\`ere une sous-vari\'et\'e $S$ connexe et totalement g\'eod\'esique de $(M,g)$, le tangent \`a $S$ est  invariant par transport par parall\'elisme le long des courbes tangentes \`a $S$. Le tangent \`a $S$ est donc en tout point du m\^eme type, $S$ est soit de type espace, soit de type temps, soit de type lumi\`ere. Par contre les feuilles d'un feuilletage totalement g\'eod\'esique n'ont aucune raison d'\^etre toutes du m\^eme type.
\subsection{Feuilletages riemanniens.}
Dans la partie \ref{mix}, il est question de flots riemanniens et de flots lorentziens. Pr\'ecisons que l'on appellera "flot" un feuilletage de dimension $1$ orient\'e, sans notion de param\'etrage. Lorsqu'un flot sera param\'etr\'e cela sera toujours explicitement pr\'ecis\'e. Nous donnons ici la d\'efinition d'un feuilletage riemannien mais le lecteur est invit\'e \`a consulter le livre de P. Molino \cite{molino} et l'article d'Y. Carri\`ere \cite {Carriere} consacr\'e aux flots riemanniens pour plus de details.
\begin{defi}
Un feuilletage $\phi$ sur une vari\'et\'e $M$ est dit (transversalement) riemannien (resp. lorentzien) s'il existe une m\'etrique riemannienne (resp. lorentzienne) $\bar \gamma$ sur $\nu(\phi)$, le fibr\'e normal de $\phi$, invariante sous l'action du pseudo-groupe d'holonomie de $\phi$. Une telle m\'etrique est dite transverse.\\
Une m\'etrique de $M$, $\gamma$, pseudo-riemannienne  sera dite quasi-fibr\'ee (par rapport \`a $\phi$) si la m\'etrique $\bar \gamma$ induite par $\gamma $ sur $\nu (\phi)$ est non d\'eg\'en\'er\'ee et invariante sous l'action du pseudo-groupe d'holonomie de $\phi$.
\end{defi}
\noindent
Nous allons aussi utiliser la notion de champ de vecteurs basique  d'un feuilletage.
\begin{defi}
Soit $\ff$ un feuilletage sur une vari\'et\'e $M$. Le champ de vecteurs $X$ nulle part tangent \`a $\ff$ est dit basique pour $\ff$ si  pour tout champ de vecteurs $Y$ tangent \`a $\ff$ on a $[X,Y]$ tangent \`a $\ff$.
\end{defi}
\subsection{Fibr\'es de Seifert.}\label{defseif}
Nous appellerons fibr\'e de Seifert (ou vari\'et\'e de Seifert) une vari\'et\'e $M$ de dimension quelconque munie d'un feuilletage riemannien de dimension $1$ transversalement orientable dont toutes les feuilles sont ferm\'ees.

La situtation est particuli\`erement agr\'eable lorsque $M$ est de dimension $3$ car les feuilles ayant de l'holonomie sont alors isol\'ees. C'est une cons\'equence directe du fait qu'une isom\'etrie directe euclidienne de $\mathbb R^2$ diff\'erente de l'identit\'e a un unique point fixe.
Ces feuilles sont appel\'ees  les fibres exceptionnelles.  La proposition II.C.3 de \cite {Carriere} donne la description d'un flot riemannien au voisinage de l'adh\'erence d'une feuille. Si toutes les feuilles sont compactes le feuilletage y est diff\'eomorphe \`a celui obtenu sur un tore solide par suspension au dessus d'un disque centr\'e en z\'ero d'une rotation $R$ d'ordre fini (d'angle rationnel) \cad au le feuilletage donn\'e par $D^2\times [0,1]/(x,0)\sim (R(x),1)$.
Le feuilletage priv\'e de ses feuilles exceptionnelles est donn\'e par une fibration localement triviale.
On retrouve la d\'efinition usuelle des fibr\'es de Seifert. Le feuilletage est vu comme une fibration (non localement triviale) au dessus de l'espace des feuilles qui n'est pas une vari\'et\'e mais un orbifold.
\hip
On  consid\'ere parfois aussi les feuilletages riemanniens non transversalement orientables,  les fibres exceptionnelles ne sont plus isol\'ees et la vari\'et\'e contient des bouteilles de Klein feuillet\'ees. Dans ce cas le flot est rev\^etu par un flot transversalement orientable.
\section{Feuilletages non d\'eg\'en\'er\'es.}\label{ndg}
On s'interesse ici aux feuilletages totalement g\'eod\'esiques dont aucune feuille n'est de type lumi\`ere. Deux cas se pr\'esentent les feuilletages dont toutes les feuilles sont de type espace ou dont toutes les feuilles sont de type temps. Ce cas admet un angle d'attaque tr\`es efficace~:
\begin{prop}\label{debut}
Soit $\ff$ un feuilletage transversalement orientable de codimension $1$.
Il existe une m\'etrique rendant $\ff$ totalement g\'eod\'esique de type espace (resp. temps) si et seulement si $\ff$ est  transverse \`a  un flot riemannien (resp. lorentzien).
\end{prop}
C'est le point de d\'epart  de l'article \cite{Car-Ghys} et de la derni\`ere partie de \cite{BMT}. Notamment on y trouve que les feuilletages g\'eod\'esiques non d\'eg\'en\'er\'es des vari\'et\'es de Seifert de dimension $3$ sont, \`a homotopie pr\`es, ceux transverses aux fibres. 
Savoir quelles sont les vari\'et\'es de Seifert munies d'un tel feuilletage \`a sucit\'e de nombreux travaux. Pour ne pas alourdir le propos nous rappelons seulement le cas o\`u $M$ est un fibr\'e en cercle, il s'agit d'un th\'eor\`eme de Milnor et Wood (cf.\cite{wood}).
\begin{theo}[Milnor-Wood]\label{milnorwood}
Soit $M$ un fibr\'e en cercle au dessus d'une surface $\Sigma$. La vari\'et\'e $M$ admet un feuilletage de codimension $1$ transverse aux fibres (\cad totalement g\'eod\'esique non d\'eg\'en\'er\'e) si et seulement si 
$$|\mathrm{eul}(M)|\leq \mathrm{max}\{0,-\chi(\Sigma)\},$$
o\`u $\mathrm{eul}(M)$ d\'esigne le nombre d'Euler de la fibration et $\chi(\Sigma)$ la caract\'eristique d'Euler de $\Sigma$.
\end{theo}
Ce r\'esultat a par la suite \'et\'e \'etendu aux fibr\'es de Seifert de dimension $3$, \`a ce sujet le lecteur pourra consulter \cite{EHN}.
\hip
Nous allons voir dans la suite ce que devient le th\'eor\`eme de Milnor et Wood, ou du moins son interpr\'etation en termes de feuilletages totalement g\'eod\'esiques,  lorsque l'on impose au feuilletage d'\^etre soit partout de type lumi\`ere, soit avec des feuilles de chaque type.
\section{Feuilletages de type lumi\`ere.}\label{lux}
\subsection{Sur les vari\'et\'es de dimension $3$.}
Les feuilletages de codimension $1$, totalement g\'eod\'esiques et de type lumi\`ere ont d\'ej\`a fait l'objet d'une \'etude par A. Zeghib \cite{Zeg-feuil}. Un tel feuilletage $\ff$ est caract\'eris\'e par le fait qu'il contient un sous-feuilletage de dimension $1$ dont la restriction \`a chaque feuille de $\ff$ est riemannienne. C'est ainsi que l'on peut voir que sur une vari\'et\'e de dimension $3$ les feuilles d'un tel feuilletage ne peuvent \^etre que des cylindres, des plans ou des tores. En particulier toute action localement libre de $\mathbb R^2$ ou de $AG$ (le groupe affine) sur une vari\'et\'e de dimension $3$ engendre un feuilletage totalement g\'eod\'esique de type lumi\`ere. 

Nous allons montrer ici que de nombreuses vari\'et\'es de Seifert ne poss\`edent pas de tels feuilletages. On commence par un r\'esultat plus g\'en\'eral sur les feuilletages dont les feuilles sont des plans, des cylindres ou des tores. On connait d\'ej\`a quelques exemples~: le feuilletage de Reeb sur $D^2\times S^1$, le feuilletage en cylindre sur $\mathcal K$ (le $I$ fibr\'e non trivial sur la bouteille de Klein), ou les feuilletages sans composantes de Reeb de $T^2\times I$ obtenus \`a partir des feuilletages sur un anneau. Ce r\'esultat ne pr\'etend pas \^etre original, il s'agit de choses connues mises ensembles.
\begin{theo}\label{gloubi}
Soit $\ff$  un feuilletage de classe $C^2$ de codimension $1$ sur une vari\'et\'e ferm\'ee de dimension $3$ dont les feuilles sont des cylindres, des plans ou des tores. Alors soit $\ff$ est minimal, soit les feuilles toriques de $\ff$ bordent des composantes diff\'eomorphes \`a $\mathcal K$, $D^2\times S^1$ ou $T^2\times I$.
\end{theo}
\pre 
Un th\'eor\`eme non publi\'e de Duminy (mais voir \cite{Duminy}) affirme que les minimaux exceptionnels d'un feuilletage de classe $C^2$ a une infinit\'e de bouts. On en d\'eduit que $\ff$ n'a pas de minimaux exceptionnels. Le feuilletage est donc soit minimal soit ses feuilles non compactes s'accumulent toutes sur ses feuilles toriques.
Dans le second cas on peut appliquer le lemme suivant issu  de \cite {rouss}, page 150~: 
\begin{lem}
Soit $\ff$ un feuilletage de codimension $1$ de classe $C^2$ d'une vari\'et\'e $M$. Soit $F$ une feuille compacte de $\ff$. Si $F$ a un voisinage ne rencontrant pas de feuilles compactes de $\ff$ et si son holonomie est ab\'elienne alors le germe de $\ff$ au voisinage de $F$ est sans holonomie en dehors de $F$.
\end{lem}
Ainsi le feuilletage $\ff$ est presque sans holonomie, c'est-\`a-dire seules les feuilles compactes peuvent avoir une holonomie non triviale. On en d\'eduit que les composantes de $\ff$ sont feuillet\'ees soit par des plans, soit par des cylindres. Ces feuilletages ont \'et\'e class\'es dans \cite {rosenberg} par R. Roussarie et H. Rosenberg dans le cas des feuilles planes et dans \cite {hector} par G. Hector dans le cas des feuilles cylindriques. $\Box$
\begin{cor}\label{feuillecompacte}
Soit $\ff$ un feuilletage lisse orientable de codimension $1$ totalement g\'eod\'esique de type lumi\`ere d'une vari\'et\'e compacte de dimension $3$. Si $\ff$ a une feuille compacte, en fait torique, alors $M$ est  un fibr\'e en tore sur le cercle ou est diff\'eomorphe \`a $\mathcal K \cup_f \mathcal K$ (vari\'et\'e obtenue en recollant 2 exemplaires de $\mathcal K$ le long de leurs bords par un diff\'eomorphisme $f$), sinon $\ff$ est minimal.
\end{cor}
\pre On sait qu'un feuilletage lisse de codimension $1$ totalement g\'eod\'esique de type lumi\`ere d'une vari\'et\'e compacte de dimension $3$ a des feuilles planes, cylindriques ou toriques. De plus il ne peut contenir de composantes de Reeb (voir \cite{Zeg-feuil}) et on suppose ici qu'il n'est pas minimal. D'apr\`es le th\'eor\`eme \ref{gloubi}, $M$ est obtenue en recollant des composantes $T^2\times I$ et $\mathcal K$, ce qui prouve le corollaire.$\Box$
\hop
Les fibr\'es en tores sur le cercle sont de rang $2$, c'est-\`a-dire qu'ils admettent une action localement libre de $\mathbb R^2$, et donc ils contiennent effectivement des feuilletages totalement g\'eod\'esiques de type lumi\`ere \`a feuilles compactes. On voudrait \'etudier maintenant les  vari\'et\'es de Seifert de dimension $3$ qui ne sont pas des fibr\'es en tores. On montre le th\'eor\`eme suivant qui concerne les fibr\'es en cercles.
\begin{theo}\label{lumos}
Si un fibr\'e en cercle $M$ orientable sur une surface orientable $\Sigma$ de genre $g$ et de nombre d'Euler $\mathrm{eul}(M)$ poss\`ede un 
feuilletage $\ff$ lisse transversalement orientable de codimension $1$ totalement g\'eod\'esique de type lumi\`ere 
alors $g\neq 0$ et si $g>1$ alors $|\mathrm{eul}(M)|= -\chi(\Sigma)=2g-2$. 
\end{theo}
\pre Supposons que $g=0$, \cad que $\Sigma=S^2$, le rev\^etement universel de $M$ est $S^3$ ou $S^2\times S^1$ et $\ff$ contient une composante de Reeb. On sait gr\^ace \`a \cite{Zeg-feuil} que c'est impossible. Lorsque $g=1$ la vari\'et\'e $M$ admet une action localement libre de $R^2$ et donc un feuilletage totalement g\'eod\'esique de type lumi\`ere. On suppose donc dor\'enavant que $g>1$.

 La suite de la preuve d\'ecoule de la simple remarque suivante. Un feuilletage de type lumi\`ere et de codimension $1$ contient un sous-feuilletage de dimension $1$ engendr\'e par l'orthogonal du tangent. La distribution d'hyperplans a donc une classe d'Euler triviale. 

D'apr\`es la proposition \ref{feuillecompacte}, $\ff$ n'a pas de feuilles compactes. Un c\'el\`ebre th\'eor\`eme de W.P. Thurston et G. Levitt \cite{thurston} nous dit qu'alors $\ff$ est homotope \`a un feuilletage transverse aux fibres.
Par ailleurs dans l'article \cite{geoghegan-nicas} R. Geoghegan et A. Nicas calculent la classe d'Euler de la distribution normale \`a un fibr\'e de Seifert "admissible"\footnote{c'est-\`a-dire orientable, de base orientable et diff\'erente de $S^2$, ou de base $S^2$ et avec au moins $3$ fibres exceptionnelles}. Son dual de Poincar\'e est \'egal \`a~: 
$$
(\chi(\Sigma) - r){\gamma_0} 
+ \sum_{j=1}^r g_{j} \in H_1(M; \mathbb Z),
$$
o\`u $\chi(\Sigma)$ repr\'esente la caract\'eristique d'Euler de la base, $r$ le nombre de fibres exceptionnelles, $\gamma_0$ la classe d'homotopie d'une fibre r\'eguli\`ere et les $g_j$ les classes des fibres singuli\`eres.

 Ici la formule se simplifie  en
$$\chi(\Sigma) \gamma_0.$$
Si la classe d'Euler est nulle  alors $\chi(\Sigma)$ est un multiple de $\mathrm{eul}(M)$, le nombre d'Euler du fibr\'e, qui est aussi l'ordre de  $\gamma'_0$. L'existence d'un feuilletage transverse entraine que $|\chi(\Sigma)|=\mathrm{eul} (M)$, cf. theoreme \ref{milnorwood}.
$\Box$
\hop
{\bf Remarque~:} Cet \'enonc\'e contient trois invariants portant le nom d'Euler. La caract\'eristique d'Euler d'une surface qu'il est inutile de pr\'esenter. Le nombre d'Euler d'un fibr\'e en cercle (ou d'un fibr\'e de Seifert) qui repr\'esente l'obstruction \`a l'existence d'une section. Et la classe d'Euler d'un fibr\'e vectoriel (ici le fibr\'e tangent au feuilletage) qui est une obstruction \`a l'existence d'une section partout non nulle \`a ce fibr\'e vectoriel. \hip
Ce r\'esultat admet-il une r\'eciproque ? Il est bien connu que lorsque $\chi(\Sigma)=\mathrm{eul} (M)$, c'est-\`a-dire lorsque $M$ est le tangent unitaire d'une surface hyperbolique, $M$ admet une action localement libre lisse du groupe affine et donc $M$ poss\`ede des feuilletages totalement g\'eod\'esiques de type lumi\`ere. Par contre si $\chi(\Sigma)=-\mathrm{eul} (M)$, aucun exemple n'est connu. On sait qu'il existe des feuilletages sans feuilles compactes sur $M$ dont les feuilles sont toutes des cylindres et des plans. Mais l'annulation de la classe d'Euler n'entraine pas l'existence d'une section. Ces feuilletages admettent-ils des sous-feuilletages ? Sont-ils g\'eod\'esiques de type lumi\`ere ?
\hip
Le th\'eor\`eme \ref{lumos} se traduit sur les fibr\'es de Seifert par un  r\'esultat \`a rev\^etement fini pr\`es.
\begin{cor}
Si $M$ est une vari\'et\'e de Seifert de dimension $3$ poss\'edant un feuilletage de classe $C^2$ totalement g\'eod\'esique de type lumi\`ere et de codimension $1$ alors $M$ a  un rev\^etement fini $\widehat M$  qui est un fibr\'e en cercles sur une surface de genre $g\geq 1$. De plus plus si $g>1$ alors $|\mathrm{eul}(\widehat M)|=2g-2$
\end{cor} 
{\bf Preuve.}  Les fibr\'es de Seifert qui ne sont pas rev\^etus par des fibr\'es en cercles n'ont pas pour rev\^etement universel $\mathbb R^3$ et donc ne poss\`edent pas de feuilletages totalement g\'eod\'esiques de type lumi\`ere. Le fibr\'e $\widehat M$  poss\`ede lui aussi un feuilletage totalement g\'eod\'esique de type lumi\`ere et le th\'eor\`eme \ref{lumos} permet de conclure.$\Box$
\hip
L'une des motivations pour \'etudier les feuilletages totalement g\'eod\'esiques de type lumi\`ere est qu'ils apparaissent lorsque l'action du groupe des diff\'eomorphismes, Diff$(M)$, sur l'espace des m\'etriques de $M$ n'est pas propre (cf. \cite{mounoud}). Cependant dans ce cas l\`a, ils ne sont que lipschitziens. La seule propri\'et\'e connue par l'auteur entrainant la pr\'esence d'un feuilletage totalement g\'eod\'esique \emph{lisse} est l'existence d'une m\'etrique analytique r\'eelle ou \`a courbure constante sur une vari\'et\'e ferm\'ee dont le groupe d'isom\'etrie est non compact (cf. \cite {zeghibII} et \cite{Zeg-feuil}). On a donc le corollaire suivant
\begin{cor}
Soit un fibr\'e en cercle $M$ orientable sur une surface $\Sigma$ de genre $g>1$ et de nombre d'Euler $\mathrm{eul}(M)$. Si $|\mathrm{eul}(M)|\neq 2g-2$ alors le groupe d'isom\'etrie de toute m\'etrique lorentzienne analytique r\'eelle de $M$ est compact.
\end{cor}
On a \'evidemment l'analogue lorsque la m\'etrique est \`a courbure constante n\'egative mais ce r\'esultat a \'et\'e montr\'e par F. Salein dans \cite{salein} en excluant aussi le cas $\mathrm{eul}(M)= 2-2g$. 
\hop
Nous allons donner rapidement un r\'esultat sur les autres vari\'et\'es de dimension $3$.
Comme on s'est retreint au cas lisse, le th\'eor\`eme 9 de \cite{Zeg-feuil} a pour corollaire~:
\begin{cor}
Soit $M$ une vari\'et\'e compacte de dimension $3$ qui n'est pas rev\^etue par un fibr\'e en tore sur le cercle. S'il existe sur $M$ un feuilletage $\ff$ lisse de codimension $1$ totalement g\'eod\'esique et de type lumi\`ere alors $\ff$ peut \^etre param\'etr\'e par une action $C^0$ du groupe affine.
\end{cor}
{\bf Preuve.} D'apr\`es le corollaire \ref{feuillecompacte} on sait que $\ff$ est minimal. \'Etant lisse on sait qu'il poss\`ede une mesure transverse si et seulement il est sans holonomie. Si $\ff$ \'etait sans holonomie  $M$ serait rev\^etu par $\mathbb T^3$ ce qui est exclu ici. Ainsi $\ff$ n'a pas de mesure transverse. Le th\'eor\`eme 9 de \cite{Zeg-feuil} affirme justement qu'un feuilletage totalement g\'eod\'esique de type lumi\`ere sans mesure transverse est param\'etr\'e par une action $C^0$ du groupe affine.$\Box$
\hip
Il semble raisonnable de penser qu'il n'existe pas de nouveaux exemples sur ces vari\'et\'es, du moins lisses.
\subsection{Un autre exemple.}
Nous donnons maintenant un r\'esultat d'existence de feuilletages totalement g\'eod\'esiques de type lumi\`ere plus \'el\'ementaire. 
\begin{prop} Il existe des feuilletages totalement g\'eod\'esiques lisses de type lumi\`ere et de codimension $1$ sur le produit de deux sph\`eres $S^{2n+1} \times S^{2m+1}$.
\end{prop}
\noindent
Il s'agit, \`a notre connaissance, du premier exemple connu de feuilletage de codimension $1$ totalement g\'eod\'esique de type lumi\`ere sur une vari\'et\'e compacte simplement connexe.\\
\pre On munit la sph\'ere $S^m$ d'un feuilletage $\ff$ de codimension $1$ et l'autre sph\`ere d'un flot riemannien $\phi$. Le feuilletage $\ff\times S^n$ est partout tangent au flot riemanien $\phi$. Ce qui d'apres \cite{Zeg-feuil} permet de conclure. $\Box$ 
\hip
{\bf Remarques~:} 1) on pourrait montrer que l'existence d'un feuilletage de codimension $1$ tangent \`a un flot riemannien sur une vari\'et\'e simplement connexe entra\^\i ne l'existence d'une action localement libre de $R^2$. Il semble donc impossible de reproduire cette construction sur les sph\`eres (cf. \cite{S1xR} sur l'existence de telles actions sur les sph\`eres). 
\hip
2) Cette proposition est \`a relier avec la question de D'Ambra et Gromov~:\\
<<Si $M$ est compacte simplement connexe l'action de Diff$(M)$ sur l'espace des m\'etriques de Lorentz est-elle propre (pour la topologie $C^2$) ?>>
\\
 Comme on l'a d\'ej\`a rappel\'e la non-propret\'e de cette action entraine l'existence d'un feuilletage lipschitzien de codimension $1$ totalement g\'eod\'esique et de type lumi\`ere (cf. \cite{mounoud} th\'eor\`eme 2.2). Ainsi sur $S^3\times S^3$ l'eventuelle pr\'esence de ces feuilletages ne permet pas de montrer cette conjecture (contrairement \`a $S^3$) mais elle ne suffit pas, bien entendu, \`a l'infirmer.
\section{Feuilletages mixtes}\label{mix}
Dans ce paragraphe nous consid\'erons les feuilletages lisses totalement g\'eod\'esiques mixtes c'est-\`a-dire poss\'edant des feuilles d'au moins deux types. Ceux ci ont \'et\'e \'etudi\'es pr\'ec\'edemment par K. Yokumoto \cite{yokumoto}. Il prouve en particulier le r\'esultat suivant.
\begin{theo}[Yokumoto] \label{theoken} Il existe des feuilletages totalement g\'eod\'esiques et de codimension $1$ sur les espaces lenticulaires. De plus si $\ff$ est une composante de Reeb d'un feuilletage totalement g\'eod\'esique alors ses feuilles planes sont de type espace et sa feuille torique est de type lumi\`ere.
\end{theo}
Nous allons donner de nouveaux exemples de tels feuilletages, principalement sur les fibr\'es de Seifert. Nous allons suivre une d\'emarche sensiblement diff\'erente.
\hip
Si on regarde le $1$-feuilletage orthogonal d'un feuilletage $\ff$ totalement g\'eod\'esique mixte, on trouve un flot riemannien sur l'ouvert des feuilles de type espace de $\ff$, lorentzien sur l'ouvert des feuilles de type temps de $\ff$ et riemannien d\'eg\'en\'er\'e sur le ferm\'e des feuilles de type lumi\`ere (cf. sections \ref{ndg} et \ref{lux}). Nous allons nous concentrer sur les situations ou les feuilles de type lumi\`ere sont isol\'ees.
Nous nous restreignons  au cas o\`u ce flot orthogonal est \emph{partout} riemannien ou lorentzien, ce qui dans ce cas est raisonnable.
Cela va nous permettre de mieux comprendre ces feuilletages, essentiellement la situation au voisinage des feuilles de type lumi\`ere et surtout d'avoir un analogue de la proposition \ref{debut}. 

Plus pr\'ecisement, nous nous posons la question : quelle position relative entre un feuilletage $\ff$ de codimension $1$ et d'un flot riemannien (ou lorentzien) $\phi$ permet de conclure qu'il existe une m\'etrique lorentzienne rendant $\ff$ totalement g\'eod\'esique et telle que $T\ff^\perp=T\phi$, en excluant bien-s\^ur le cas non d\'eg\'en\'er\'e o\`u $\phi$ et $\ff$ sont transverses. On trouve dans l'article de K. Yokumoto \cite{yokumoto}  des exemples de feuilletages totalement g\'eod\'esiques invariants par une action de $S^1$. On cherche \`a relacher cette condition au maximum.

L'hypoth\`ese minimale (et qui semble raisonnable) est de supposer que pour chaque feuille $F$ de $\ff$ soit $\phi$ est tangent \`a $F$ partout soit $\phi$ est transverse \`a $F$ partout. 
Mais nous allons voir que ce ne sera pas suffisant. Posons donc une premi\`ere d\'efinition.
\begin{defi}
Soit $\ff$ un feuilletage de codimension $1$ et $\phi$ un feuilletage de dimension $1$. On dira que $\ff$ et $\phi$ sont transverses ou tangents feuille \`a feuille si pour toute feuille $F$ de $\ff$, $F$ et $\phi$ sont soit partout transverses, soit partout tangents.

On dira qu'un feuille $F$ de $\ff$ est une feuille de tangence si en tout point $x$ de $F$ on a $T_x\phi\subset T_x\ff$.
\end{defi}
Commen\c cons par une remarque sur les flots lorentziens. Supposons que $\phi$ soit un flot lorentzien tangent \`a une hypersurface $S$ de type espace pour une m\'etrique quasi-fibr\'ee. Le flot $\phi$ pr\'eserve la direction orthogonale \`a cette hypersurface qui est de type temps et donc en tout point de $S$ sa differentielle est born\'ee. Il d\'ecoule du th\'eor\'eme B de \cite{BMT} que la diff\'erentielle de $\phi$ est partout born\'ee et que $\phi$ est riemannien. Cette situation est exactement celle rencontr\'ee le long d'une feuille de tangence (voir proposition \ref{globi}). C'est pourquoi nous n'\'etudierons que le cas des feuilletages transverses ou tangents \`a des flots riemanniens. L'hypoth\`ese lorentzien sera toujours en plus. Il est sans doute bon de rappeler qu'un flot \`a la fois lorentzien et riemannien est un flot riemannien poss\'edant un champ de vecteurs basique obtenu en diagonalisant simultan\'ement les m\'etriques transverses.
\subsection{Un premier obstacle.}\label{prems}
Regardons une premi\`ere propri\'et\'e globale des feuilletages totalement g\'eod\'esiques. Pour cela nous donnons une nouvelle d\'efinition. Un feuilletage de codimension $1$ dont l'orthogonal est orient\'e est transversalement orient\'e (m\^eme si certaines feuilles sont de type lumi\`ere), on supposera donc  toujours $\ff$  transversalement orientable dans la suite.
\begin{defi}
Soit $\ff$ un feuilletage de codimension $1$ transversalement orientable et transverse ou tangent feuille \`a feuille \`a  un flot $\phi$. Soit  $F_0$ une feuille de tangence de $\ff$ et $X$ une trivialisation de $T\phi$. \`A
 tout champ de vecteurs $Z$ transverse \`a $\ff$ on peut associer un champ de vecteurs $Y$ non nul et appartenant \`a $T\ff\cap \mathrm {Vect}(X, Z)$. Il existe alors 2 fonctions $a$ et $b$ telle que $Y=a X + bZ$. On dira que $F_0$ est \emph{attractive vue de $\phi$} si $b$ change de signe au voisinage de $F_0$. 
\end{defi}
\noindent Les champs de vecteurs non nuls transverses \`a $\ff$ sont tous homotopes et donc cette propri\'et\'e ne d\'epend pas du choix de $Z$. La proposition qui suit montre comment cette propri\'et\'e conditionne le type des feuilles des feuilletages totalement g\'eod\'esiques.
\begin{prop}\label{attractive}
Soit $\ff$ un feuilletage totalement g\'eod\'esique sur une vari\'et\'e lorentzienne $(M,g)$. Supposons que $F_0$ soit une feuille de type lumi\`ere isol\'ee (ces voisines sont non d\'eg\'en\'er\'ees). Alors les feuilles voisines de $F_0$ sont toutes de m\^eme type si et seulement si $F_0$ n'est pas attractive vue de $\ff^\perp$.
\end{prop}
\pre Soit $X$ une trivialisation de $\ff^\perp$ et soit $Z$ un champ de vecteurs de type temps transverse \`a $\ff$ et tel que $g(Z,Z)=-1$. Comme ci-dessus on consid\'ere le champ $Y$ tel que $Y \subset T\ff$ et  que $Y=X+b Z$, $b$ \'etant une fonction lisse. On a $g(Y,Y)=b-b^2$ dont le signe ne d\'epend que de l'attractivit\'e vue de $\ff^\perp$.$\Box$
\hop

Cette proposition nous permet d\'ej\`a de donner un contre-exemple \`a notre hypoth\`ese "minimale". En effet consid\'erons un feuilletage de Reeb de $S^3$ (\cad un feuilletage constitu\'e de deux composantes de Reeb) invariant sous l'action du champ de Hopf, qui d\'efinit un flot riemannien $\phi$, et  transverse ou tangent feuille \`a feuille avec celui-ci. \`A l'int\'erieur de chacun des tores solides on a la possibilit\'e de retourner la composante de Reeb car l'holonomie de la feuille torique est $C^\infty$ infinit\'esimalement triviale (\cad le developpement de Taylor de son holonomie est trivial \`a tout ordre). On peut donc supposer que la feuille torique est attractive vue de $\phi$. Le th\'eor\`eme \ref{theoken} avec la proposition \ref{attractive} nous permet d'affirmer qu'il n'existe pas de m\'etriques lorentziennes rendant ce feuilletage de Reeb totalement g\'eod\'esique.
\hip
De m\^eme si $\ff$ poss\`ede une unique feuille de tangence avec $\phi$, si elle ne s\'epare pas la vari\'et\'e et si elle est  attractive vue de $\phi$ alors $\ff$ ne peut pas \^etre rendue totalement g\'eod\'esique.
On voit qu'il est n\'ecessaire de rajouter des hypoth\`eses simplement d'ordre combinatoire (on doit avoir un nombre pair de changement de signe pour pouvoir recoller) et d'autres plus subtiles sur les propri\'et\'es de $\phi$ qui doit \^etre lorentzien sur certaines composantes de $\ff$.
Par exemple, pour n'avoir que des feuilles de type espace ou de type lumi\`ere il est n\'ecessaire d'imposer \`a  $\ff$ de n'avoir pas de feuilles de tangences avec $\phi$ qui soit attractives vues de $\phi$.
\hip
Apr\`es ces obstructions de type global, nous allons voir que localement, \cad au voisinage d'une feuille de tangence, on ne peut pas non plus toujours rendre le feuilletage g\'eod\'esique.
\subsection{Au voisinage d'une feuille de tangence.}\label{local}
\'Etant donn\'e un feuilletage $\ff$ de codimension $1$ et un flot riemannien $\phi$ on cherche \`a construire une m\'etrique $g$ rendant $\ff$ totalement g\'eod\'esique. 
Notre connaissance des cas non-d\'eg\'en\'er\'es et de type lumi\`ere nous permet de donner imm\'ediatement le crit\'ere suivant. Crit\`ere que nous utiliserons par la suite.
On retrouve ce r\'esultat,  formul\'e diff\'eremment dans  \cite{yokumoto} proposition 2.13.
\begin{prop} \label{critere}
Soit $\ff$ un feuilletage de codimension $1$ transversalement orientable $\phi$ un flot. On d\'enote $U_T$ l'ouvert de $M$ sur lequel ces feuilletages sont transverses.
Il existe une m\'etrique pseudo-riemannienne $g$ telle que $\ff$ et $\phi$ sont orthogonaux qui rend $\ff$ totalement g\'eod\'esique si et seulement si $\ff$ et $\phi$ sont transverses ou tangents feuille \`a feuille, la m\'etrique $g$ est quasi-fibr\'ee par rapport \`a $\phi$ sur $U_T$ et la restriction de $g$ aux feuilles de tangences est elle aussi quasi-fibr\'ee par rapport \`a $\phi$.
\end{prop}
Tous le probl\`eme dans la suite sera de reussir \`a construire une telle m\'etrique.
Nous allons voir que m\^eme localement on rencontre des difficult\'es. Il va nous faloir d\'ecrire la position relative de $\ff$ et $\phi$, pour cela nous avons tout d'abord besoin du lemme suivant.
\begin{lemme}\label{Z}
 Soit $F_0$ une hypersurface plong\'ee tangente \`a un flot riemannien $\phi$. Il existe au voisinage de $F_0$ un champ de vecteurs $Z$ transverse \`a $F_0$ et pr\'eserv\'e par $\phi$ (basique).
\end{lemme}
\pre Il suffit de consid\'erer pour une  m\'etrique quasi-fibr\'ee, l'orthogonal de $F_0$ et les g\'eod\'esiques issues de cet orthogonal.$\Box$
\hip
On peut maintenant formuler la d\'efinition suivante.
\begin{defi}\label{adapte}
Soit $\ff$ un feuilletage de codimension $1$ transversalement orientable et transverse ou tangent feuille \`a feuille \`a un flot riemannien $\phi$. Soit $F_0$ une feuille de tangence propre (plong\'ee) et isol\'ee, $Z$ le champ de vecteurs fourni par le lemme \ref{Z}, $X$ une trivialisation de $T\phi$. Soit $Y$ un champ de vecteurs  tangent \`a $T\ff\cap \mathrm{Vect}(X,Z)$, il existe deux fonctions $a$ et $b$ telles que $Y=aX+bZ$ (on peut m\^eme supposer que $a^2+b^2=1$).

On dira que $\ff$ est \emph{adapt\'e \`a $\phi$ en $F_0$} s'il existe, au voisinage de $F_0$, une fonction $\beta$ telle que $X.\beta=0$ et que $b/\beta$ soit lisse et ne s'annulle pas.

On dira que $\ff$ est adapt\'e \`a $\phi$ s'il l'est pour toute feuille de tangence.
\end{defi}
Cette condition  ne d\'epend que des germes du feuilletage $\ff$ en les feuilles de tangence.
Cette fonction $b$ d\'ecrit la fa\c con dont le feuilletage $\ff$ se comporte au voisinage de la feuille $F_0$ dans une direction du moins. Il serait int\'eressant de donner une interpr\'etation plus g\'eom\'etrique, en relation avec l'holonomie du feuilletage par exemple. Ceci est faisable lorsque les feuilles de $\phi_{|_{F_0}}$ sont toutes ferm\'ees mais cela semble difficile dans le cas g\'en\'eral.\\
Clairement lorsque $b$ est constante le long de $\phi$ cette hypoth\`ese est v\'erifi\'ee. C'est le cas par exemple si on peut param\'etrer $\phi$ en $\phi_t$ de fa\c con isom\'etrique et si $\ff$ est invariant par $\phi_t$. \`A partir d'une telle situation, on peut perturber $\ff$ loin des feuilles de tangences sans modifier le fait que $\ff$ est adapt\'e \`a $\phi$. Nous verrons donc la propri\'et\'e <<\^etre adapt\'e \`a $\phi$>> comme les d\'eformations maximales autoris\'ees depuis une situation invariante ou comme une propri\'et\'e  faible d'invariance de $\ff$. Notons que cette propri\'et\'e n'est pas triviale, on construit facilement des feuilletages transverses ou tangents feuille \`a feuille \`a un flot riemannien mais non adapt\'es \`a celui-ci.
\hip
Cette d\'efinition est toutefois le bon crit\`ere comme le montre la proposition suivante~:
\begin{prop}\label{voisinage}
Soit $\ff$ un feuilletage de codimension $1$ transversalement orientable et transverse ou tangent feuille \`a feuille \`a un flot riemannien $\phi$. Soit $F_0$ une feuille de tangence compacte et isol\'ee. Il existe au voisinage de $F_0$ une m\'etrique lorentzienne $g$  telle que  le tangent \`a $\phi$ soit l'orthogonal du tangent \`a $\ff$ rendant $\ff$ totalement g\'eod\'esique si et seulement $\ff$ est adapt\'e \`a $\phi$ en $F_0$.
\end{prop}
\pre On va construire, au voisinage de $F_0$, une m\'etrique lorentzienne quasi-fibr\'ee sur $\{x\in M\, |\, T_x\phi + T_x \ff=T_x M\}$ et dont la restriction \`a $F_0$ sera encore quasi-fibr\'ee. Cette m\'etrique rendra effectivement le feuilletage $\ff$ totalement g\'eod\'esique.
 Par hypoth\`ese $F_0$ est compacte et il existe un voisinage de $F_0$ sur lequel $F_0$ est la seule feuille de tangence entre $\ff$ et $\phi$. D'apr\`es le lemme \ref{Z}, il existe un champ de vecteurs $Z$ invariant par $\phi$ et transverse \`a $F_0$. Soit $\gamma$ une m\'etrique riemannienne quasi-fibr\'ee par rapport \`a $\phi$. Soit $X$ une trivialisation du tangent \`a $\phi$.
Consid\'erons le champ de vecteurs $Y$ introduit dans la d\'efinition \ref{adapte}, on a $Y=aX+bZ$ et $a^2+b^2=1$.
Au voisinage de $F_0$, on va d\'efinir une m\'etrique lorentzienne $g$. On note $P$ la distribution  orthogonale, pour $\gamma$, \`a  $\mathrm{Vect}(X,Z)$.
Consid\'erons la d\'ecomposition du tangent \`a la vari\'et\'e suivante : $TM=\mathbb R X \oplus \mathbb R Z\oplus P$. Prenons la m\'etrique qui se repr\'esente dans une base adapt\'ee \`a cette d\'ecomposition par la matrice suivante~:
$$\begin {pmatrix}
-\frac{b^2}{\beta}  & \frac{a\,b}{\beta}& 0\\
\frac{a\,b}{\beta}   & \frac{b^2}{\beta}   & 0\\
   0          &        0   & \gamma_{|_P}
\end{pmatrix},
$$
o\`u $\beta$ est la fonction introduite lors de la d\'efinition \ref{adapte}. Cette m\'etrique est bien lorentzienne et lisse de plus sa restriction \`a  $TF_0$ a bien les propri\'et\'es voulues.
En un point $x$ tel que $b(x)\neq 0$, c'est-\`a-dire en un point n'appartenant pas \`a $F_0$, on peut d\'ecomposer le tangent \`a la vari\'et\'e en $\mathbb R X\oplus\mathbb R Y \oplus P$ (de plus $T\ff=\mathbb R Y \oplus P$) la m\'etrique s'\'ecrit alors~:
$$\begin {pmatrix}
-\frac{b^2}{\beta}  & 0& 0\\
0   & \frac{b^2}{\beta}   & 0\\
   0          &        0   & \gamma_{|_P}
\end{pmatrix}.
$$
Clairement $T\phi^\perp=T\ff$. Le flot riemannien $\phi$ respecte la d\'ecomposition orthogonale $\mathbb RY\oplus P$ et donc $g$ est quasi-fibr\'ee si et seulement si  le champ de vecteurs $Y/\sqrt{|g(Y,Y)|}$ est pr\'eserv\'e par $ \phi$, \cad est basique. Ce qui est vrai si et seulement si $\beta$ est constante le long de $\phi$ ce que l'on a suppos\'e vrai. On a donc une m\'etrique $g$ d\'efinie  au voisinage de $F_0$  telle que $\ff$ est totalement g\'eod\'esique et $T\phi^\perp=T\ff$.
\hip
R\'eciproquement supposons qu'il existe une m\'etrique lorentzienne $g$ telle que $T\ff$ soit orthogonal \`a $T\phi$ et que $\ff$ soit totalement g\'eod\'esique. On consid\`ere au voisinage de la  feuille de tangence $F_0$ les champs de vecteurs $X,\,Y$ et $Z$ comme pr\'ec\'edemment. Si on est suffisament proche de $F_0$, $Z$ n'est pas tangent \`a $\ff$  et donc le champ de plans Vect$(X,Z)$ n'est pas d\'eg\'en\'er\'e. Soit $P$ la distribution de plans de codimension $2$ orthogonale, pour $g$, \`a Vect$(X,Z)$, $P$ est bien suppl\'ementaire de Vect$(X,Z)$.
 Comme plus haut on \'ecrit $Y=aX+bZ$ avec $a^2+b^2=1$ et on peut repr\'esenter $g$, dans une base adapt\'ee \`a la d\'ecomposition du tangent en $\mathbb R X \oplus \mathbb R Z \oplus P$, par la matrice~:
$$\begin{pmatrix}
-b\, \alpha & a\,\alpha&0\\
a\,\alpha   & c\alpha & 0       \\
0&0&g_{|_P}
\end{pmatrix},$$
o\`u $\alpha$ et $c$ sont deux fonctions. On a 
$g(Y,Y)= a^2\alpha b + c\, b^2 \alpha$. La fonction $c$ \'etant born\'ee et la fonction $\alpha$ partout non nulle  on voit que $Y$ ne peut pas \^etre de type lumi\`ere (sauf sur $F_0$). Le feuilletage $\ff$ \'etant totalement g\'eod\'esique et  $Z$ \'etant  pr\'eserv\'e par $\phi$, on en d\'eduit que le champ de vecteurs $Y/\sqrt{|g(Y,Y)|}$, d\'efini si $b\neq 0$, \cad hors de $F_0$, est lui aussi invariant par $\phi$ et donc $X.\big(g(Y,Y)/b^2\big)=0$. On a donc 
$$X.\big(\alpha (1/b-b+c)\big)=0.$$
Sachant que $\alpha$ ne s'annule pas et que $c$ est born\'ee on en d\'eduit que la fonction $\beta=\left(\alpha(1/b-b+c)\right)^{-1}$ v\'erifie  $X.\beta=0$ et   $b/\beta$ est lisse et ne s'annule pas sur $F_0$. Le feuilletage $\ff$ est bien adapt\'e \`a $\phi$ en $F_0$.
$\Box$
\hop
L'exemple construit au paragraphe \ref{prems} est bien adapt\'e \`a son orthogonal mais cela ne suffit pas.
\subsection{O\`u l'on r\'esoud le probl\`eme global.}
Il nous reste \`a r\'esoudre les difficult\'es rencontr\'ees au paragraphe \ref{prems}. C'est ce que nous permet de faire la proposition suivante.
\begin{prop}\label{globi}
Soit $\ff$ un feuilletage totalement g\'eod\'esique de codimension $1$, sur une vari\'et\'e lorentzienne $(M,g)$, dont l'orthogonal est un flot riemannien $\phi$ et dont les feuilles de type lumi\`ere sont compactes et en nombre fini non nul. 
Soit $K$ une composante de $\ff$ de type temps maximale \cad telle que tout voisinage du bord de $K$ rencontre  des feuilles de type espace. 
Il existe alors un voisinage $U$ de $K$ tel que $\phi$ restreint \`a $U$ soit (transversalement)  lorentzien. De plus il existe sur $U$ une m\'etrique quasi-fibr\'ee $h$ transversalement lorentzienne telle que les feuilles  de $\ff$ qui sont de type lumi\`ere pour $g$ sont de type espace pour $h$.\end{prop}
\pre
Nous allons prouver qu'il existe un voisinage $U$ de $K$ sur lequel $\phi$ pr\'eserve une direction, cela montrera bien que $\phi$ y est lorentzien.
 Il est clair que $\phi$ restreint \`a l'int\'erieur de $K$ est lorentzien (m\^eme si l'int\'erieur de $K$ contient des feuilles de type lumi\`ere) : la direction obtenue en diagonalisant simultan\'ement $g$ restreinte au tangent \`a $\ff$ et une m\'etrique  riemannienne transverse donne la direction invariante d\'esir\'ee. 
D'autre part le bord de $K$ est compos\'e de feuilles de type lumi\`ere, le lemme \ref{Z} donne une autre direction invariante sur ce voisinage.

Reprenons les notations de la deuxi\`eme partie de la preuve de la proposition \ref{voisinage}. 
La m\'etrique $g_{|_P}$ est forc\'ement riemannienne (sinon la m\'etrique n'est plus lorentzienne lorsque $b=0$). Ainsi quitte \`a modifier la m\'etrique riemannienne transverse on peut supposer que la direction transverse donn\'ee par $Z$  est bien la m\^eme que celle obtenue plus haut. On a bien une m\'etrique lorenzienne transverse dont la restriction \`a $P$ est riemannienne.$\Box$
\hip
Cette proposition nous invite \`a d\'efinir la notion suivante~:
\begin{defi}\label{C}
Soit $\ff$ un feuilletage de codimension $1$  et $\phi$ un flot riemannien  transverses ou tangents 
feuille \`a feuille et n'ayant qu'un nombre fini de feuilles de tangence. On dira que $\ff$ et $\phi$ sont \emph{globalement compatibles} s'il existe un ouvert $U$ de $M$ dont la fronti\`ere ne rencontre pas les feuilles de tangences et une m\'etrique quasi-fibr\'ee transversalement lorentzienne $h$ sur $U$ telle que les feuilles de tangences de $\ff$ contenues dans $U$ sont de type espace pour $h$, que tout chemin traversant une composante 
connexe de $U$ rencontre $2$ feuilles attractives vues de $\phi$ et que toutes les feuilles de 
tangences attractives vues de $\phi$ sont dans $U$. 
\end{defi}
On voit bien que malgr\`es son apparente opacit\'e cette d\'efinition est incontournable et n'est probablement pas simplifiable. Mise avec les r\'esultats du paragraphe \ref{local}, elle nous permet de formuler ce th\'eor\`eme~:
\begin{theo}\label{lapinot}
Soit $\ff$ un feuilletage de codimension $1$ transversalement orientable et transverse ou tangent feuille \`a feuille \`a un flot riemannien $\phi$ sur une vari\'et\'e $M$. On suppose de plus que les feuilles de tangences de $\ff$ sont compactes et en nombre fini.
Il existe une m\'etrique lorentzienne $g$ rendant $\ff$ totalement g\'eod\'esique et telle que  $T\ff$ soit orthogonal \`a $T\phi$ si et seulement si $\ff$ est  adapt\'e \`a $\phi$ et si $\ff$ et $\phi$ sont globalement compatibles.
\end{theo}
\pre 
Tout d'abord il est clair que sous les hypoth\`eses du th\'eor\`eme, si $\ff$ est totalement g\'eod\'esique avec $T\ff$  orthogonal \`a $T\phi$ alors $\ff$ est  adapt\'e \`a $\phi$ et $\ff$ et $\phi$ sont globalement compatibles : il s'agit des propositions \ref{voisinage} et \ref{globi}. \'Etudions la r\'eciproque.
\hip
On suppose que $\ff$ est  adapt\'e \`a $\phi$ et que $\ff$ et $\phi$ sont globalement compatibles. On reprendra certaines notations de la preuve de la proposition \ref{voisinage}.
D'apr\`es la proposition \ref{voisinage}, il existe une m\'etrique $g$ ayant les propri\'et\'es voulues mais d\'efinie au voisinage de chacune des feuilles de tangences seulement.
 Il reste \`a recoller les morceaux sans perdre ces propri\'et\'es.
Consid\'erons $U$ une composante connexe de $M$ priv\'ee des feuilles de tangences de $\ff$. Le fait que $\ff$ et $\phi$ soient globalement compatibles nous permet de choisir $g$ de telle sorte qu'au voisinage du bord de $U$, la restriction de  $g$ \`a  $T\ff$ est soit riemannienne,  soit lorentzienne. 
Les hypoth\`eses faites sur le nombre de feuilles attractives vues de $\phi$, reviennent \`a  supposer que si la restriction de $g$ \`a $\ff$ est lorentzienne alors $\phi$ est transversalement lorentzienne sur un voisinage $V$ de $U$, ce qui est n\'ecessaire au raccord de $g$ d'apr\`es la proposition \ref{globi}. Le champ de vecteurs $Z$ qui a servi \`a la construction de $g$ au voisinage des feuilles de tangence peut alors \^etre \'etendu \`a tout $V$ (mais pas la fonction $\beta$ a priori). 

Sur l'ouvert $U$ le flot $\phi$ et le feuilletage $\ff$ sont transverses, on peut donc consid\'erer une m\'etrique $g'$ quasi-fibr\'ee telle que $T\phi$ soit orthogonal \`a $T\ff$. On ne cherche pas \`a \'etendre $g'$ au bord de $U$. Si la restriction de $g$ \`a $\ff$ est lorentzienne sur $U$, on choisit $g'$ de telle sorte qu'elle co\"\i ncide avec la m\'etrique $h$ donn\'ee par la proposition \ref{globi} sur le fibr\'e normal de $\phi$ d\'enot\'e $\nu(\phi)$ (on a $\nu(\phi)=TM/T\phi$). Si $g$ restreint \`a $\ff$ est riemannienne, on prend une m\'etrique $g'$ dont la restriction \`a $\nu(\phi)$ est riemannienne.
Dans les deux cas, si $\zeta$ est une fonction plateau s'annulant hors du domaine de d\'efinition de $g$ et valant $1$ sur un voisinage du bord de $U$ alors la m\'etrique $g_{\zeta}=\zeta g+ (1-\zeta)g'$ est une m\'etrique lorentzienne lisse pour laquelle  $T\phi$ est bien orthogonal \`a $T\ff$. 
De plus il est clair que $g_\zeta$ est quasi-fibr\'ee si et seulement si
$\zeta$ est constante le long de $\phi$.
\hip
Les deux situations se ram\`ene donc au m\^eme probl\`eme : trouver une fonction plateau v\'erifiant les conditions ci-dessus. Pour cela on revient dans la vari\'et\'e de d\'epart $M$ munie d'une m\'etrique \emph{riemannienne} quasi-fibr\'ee $\gamma$. Au voisinage d'une feuille de tangence $F_0$ on consid\`ere la fonction $\delta$ qui a un point $x$ de $M$ associe la distance de $x$ \`a $F_0$. Comme $F_0$ est feuillet\'ee par $\phi$ et que $\gamma$ est quasi-fibr\'ee cette fonction est constante le long de $\phi$. On construit alors ais\'ement la fonction $\zeta$ \`a partir de $\delta$.
$\Box$\hip
On peut donner le corollaire suivant, qui est plut\^ot une relecture de la preuve pr\'ec\'edente sous des hypoth\`eses plus fortes.
\begin{cor}\label{invariant}
Sous les  hypoth\`eses du th\'eor\`eme \ref{lapinot}. Supposons de plus que $\phi$ peut \^etre param\'etr\'e en un champ de Killing riemannien $\phi_t$ et que $\ff$ est invariant par $\phi_t$.
Alors il existe une m\'etrique lorentzienne $g$ \emph{invariante} par $\phi_t$ rendant $\ff$ totalement g\'eod\'esique si et seulement si $\ff$ et $\phi$ sont globalement compatibles.
\end{cor}
\pre Sous ces conditions, on peut prendre $\beta=b$ et donc $\ff$ et adapt\'e \`a $\phi$. Le th\'eor\`eme \ref{lapinot} s'applique. En rep\'etant la preuve on voit bien que les diff\'erents ingr\'edients de la preuve sont tous invariants par $\phi_t$ : les fonctions $a$ et $b$, la distribution $P$ et $\gamma_{|_P}$.  Le flot param\'etr\'e $\phi_t$ pr\'eserve donc la m\'etrique ainsi construite.
$\Box$
\subsection{O\`u l'on rencontre enfin des feuilletages totalement g\'eod\'esiques.}
Nous allons donner un exemple illustrant la propri\'et\'e \ref{C} de compatibilit\'e globale. Consid\'erons le  feuilletage de Reeb sur $S^3$ vu plus haut au paragraphe \ref{prems} , celui qui n'est pas g\'eod\'esible. Prenons comme flot riemannien $\phi$, param\'etr\'e en $\phi_t$, le champ de Hopf. Au centre de l'une de ses composantes de Reeb on peut effectuer un tourbillonement de Reeb (cf. \cite{Ca-Co} p. 89). On obtient un nouveau feuilletage $\ff$ invariant par $\phi_t$, compos\'e de 2 composantes de Reeb et d'un feuilletage par plans de $T^2\times I$. Les 2 feuilles toriques sont attractives vues de $\phi$ et, restreint \`a un voisinage de $T^2\times I$, le flot $\phi$ est bien lorentzien (\cad transversalement parall\'elisable en dimension $3$).
 
Ce feuilletage satisfait donc aux hypoth\`eses du corollaire \ref{invariant}, il existe donc une m\'etrique $g$ rendant $\ff$ totalement g\'eod\'esique. Les feuilles des composantes de Reeb sont alors de type espace, celles de la composante $T^2\times I$  de type temps et les 2 feuilles compactes de type lumi\`ere.

Voyons comment ce m\^eme genre de construction permet de donner de nombreux exem\-ples. Nous donnons un premier 
r\'esultat qui illustre comment cr\'eer de nouveaux exemples \`a partir d'anciens.
\begin{prop}
Soit $M$ une vari\'et\'e de Seifert de dimension quelconque.
Si $M$ poss\`ede un feuilletage transverse aux fibres (donc totalement g\'eod\'esible non d\'eg\'en\'er\'e)  alors il poss\`ede aussi une infinit\'e de feuilletages totalement g\'eod\'esiques mixtes de codimension $1$ deux \`a deux non diff\'eomorphes.
\end{prop}
\pre
Les feuilletages sont simplement obtenus en effectuant des tourbillonements de Reeb au voisinage d'autant de fibres r\'eguli\`eres (de feuilles sans holonomie) que l'on souhaite. Quitte \`a retourner ces composantes de Reeb, on voit que les feuilletages obtenus v\'erifient bien les hypoth\`eses du th\'eor\`eme \ref{lapinot}.
$\Box$
\hip
Nous allons maintenant donner des exemples sur des vari\'et\'es ne poss\'edant pas de feuilletages g\'eod\'esiques non d\'eg\'en\'er\'es.
\begin{theo}\label{seifert}
Quitte \`a prendre un rev\^etement \`a $2$ feuillets, tout fibr\'e de Seifert de dimension $3$,
 admet un feuilletage totalement g\'eod\'esique.
\end{theo}
\pre Soit $M\rightarrow B$ un fibr\'e de Seifert. On a vu au paragraphe \ref{defseif} que cette fibration provient d'un feuilletage riemannien de dimension $1$. Si n\'ecessaire nous prenons un rev\^etement \`a deux feuillets et nous supposons que l'on est bien en pr\'esence d'un \emph{flot} riemannien \cad que ce feuilletage est bien orient\'e. On commence alors la construction d'un feuilletage $\ff$.
Il existe une collection finie de boules de $B$ ferm\'ees disjointes, d\'enot\'ees  $D_i$,   telle que le fibr\'e restreint \`a $B\setminus \cup_i D_i$ soit localement trivial et admette une section, c'est-\`a-dire soit trivial.  Le fibr\'e obtenu admet donc un feuilletage  transverse au bord et invariant sous l'action de $S^1$. Les feuilletages induits sur les composantes connexes du bord de $(B\setminus \cup_i D_i)\times S^1$ sont triviaux on peut donc modifier le feuilletage de telle sorte qu'il devienne tangent au bord (op\'eration appel\'ee "spinning" dans \cite {Ca-Co}, p.84) et reste invariant sous l'action de $S^1$. Sur chaque $D_i\times S^1$, on met un feuilletage de Reeb  invariant par translation le long de $S^1$ et par rotation sur $D_i$ et dont l'holonomie de la feuille compacte  est infinit\'esimalement $C^\infty$ triviale. Pour obtenir un feuilletage de $M$, il suffit de recoller ces composantes de Reeb le long des composantes du bord. Ceci est possible gr\^ace \`a l'hypoth\`ese sur l'holonomie des feuilles compactes, cette hypoth\`ese  permet aussi de retourner, si n\'ecessaire, les composantes de Reeb et de s'assurer ainsi que les feuilles compactes ne sont pas attractives vues des fibres de $M$ (on pourrait aussi faire un tourbillonement de Reeb au centre de la composante de Reeb comme pr\'ec\'edemment). Le feuilletage satisfait toujours \`a l'hypoth\`ese d'invariance et donc aux hypoth\`eses du th\'eor\`eme \ref{lapinot}. Il existe donc un feuilletage totalement g\'eod\'esique sur $M$.$\Box$
\hip
Cette preuve s'applique sans modifications aux fibr\'es de Seifert de dimension sup\'erieure ayant des fibres exceptionnelles isol\'ees donc en particulier s'il n'y a pas de fibres exceptionnelles~:
\begin{theo}\label{cercle}
Si $M$ est un fibr\'e en cercle  orientable alors $M$ poss\`ede un feuilletage totalement g\'eod\'esique mixte de codimension $1$.
\end{theo}
Ce th\'eor\`eme s'applique par exemple aux sph\`eres de dimension impaire. Elles poss\`edent donc toutes des feuilletages totalement g\'eod\'esiques.
\hop
La vari\'et\'e $T^3_A=T^2\times [0,1]/\left((x,0)\sim (Ax,1)\right)$ o\`u $A$ un automorphisme hyperbolique du tore poss\`ede un flot riemannien bien qu'elle ne soit pas un fibr\'e de Seifert. On sait d'apr\'es \cite{Car-Ghys} et \cite{BMT} que $T_A^3$ poss\`ede des feuilletages de codimension $1$ totalement g\'eod\'esiques non d\'eg\'en\'er\'es (de type espace et de type temps) et par \cite{Zeg-feuil} qu'il en poss\'ede de type lumi\`ere ($T^3_A$ est de rang $2$). Montrons qu'il en poss\`ede aussi de mixtes.
 
Consid\'erons un feuilletage $\ff$ de $T^3_A$ sans composantes de Reeb et avec un nombre fini de feuilles compactes,  on suppose de plus que ces feuilles compactes sont des fibres de la fibration en tores sur le cercle et que les feuilles non compactes sont transverses \`a ces fibres.
Consid\'erons $\phi$ le flot riemannien sur $T^3_A$ tangent aux fibres, $\ff$ et $\phi$ sont tangents ou transverses feuille \`a feuille. Au voisinage d'un feuille compacte, on consid\`ere une trivialisation de la fibration. On est donc sur $T^2\times ]0,1[$ et $\phi$ est un feuilletage sur les tores de pente constante irrationnelle. On peut modifier $\ff$ de telle sorte que dans cette trivialisation il soit de plus invariant par la translation correspondante. Le feuilletage $\ff$ est alors adapt\'e \`a $\phi$, de plus  $\phi$ est transversalement parall\'elisable et donc riemannien et lorentzien. Si de plus $\ff$ a un nombre pair de feuilles attractives vues de $\phi$ alors le couple $(\ff,\, \phi)$ satisfait bien aux hypoth\`eses du th\'eor\`eme \ref{lapinot}.

On peut facilement g\'en\'eraliser cette construction aux vari\'et\'es de dimension sup\'erieure poss\'edant un flot riemannien dont l'adh\'erence de chaque feuille est un tore de codimension~$1$. On montre :
\begin{prop}
Soit $M$  une vari\'et\'e compacte de dimension $n\geq 3$ poss\'edant un flot riemannien transversalement orientable  $\phi$. Si $\phi$ poss\`ede une feuille dont l'adh\'erence est de codimension $1$ alors $M$ poss\`ede un feuilletage totalement g\'eod\'esique mixte de codimension~$1$.
\end{prop}
\pre
Deux cas se pr\'esentent. Soit l'adh\'erence de chaque feuille est un tore de codimension $1$ et $M$ est un fibr\'e en tores sur le cercle et on reproduit la construction ci-dessus. 
Soit il existe une feuille dont l'adh\'erence est de codimension plus grande. Cette situation est d\'ecrite dans \cite{molino}. On y voit que dans ce cas  les adh\'erences des feuilles sont des tores $T^{n-1}$ sauf deux d'entre elles qui sont des tores $T^{n-2}$. Le long de chacun de ces tores le flot $\phi$ est conjugu\'e \`a un flot lin\'eaire dense.
La vari\'et\'e $M$ est obtenue en recollant deux exemplaires de $D^2\times T^{n-2}$ le long de leur bord. 
On d\'ecompose $T^{n-2}$ en $S^1\times T^{n-3}$ et on prend $\ff_0=\mathcal R\times T^{n-3}$, o\`u $\mathcal R$ est un feuilletage de Reeb de $D^2\times S^1$ transverse au facteur $S^1$, invariant sous son action et tel que l'holonomie de la feuille compacte soit $C^\infty$ infinit\'esimalement triviale. Le flot $\phi$ est tangent \`a $\ff$ sur le bord de $D^2\times T^{n-2}$ et transverse sur l'int\'erieur. On recolle deux exemplaires de ce feuilletage, on obtient  un feuilletage $\ff$ adapt\'e \`a $\phi$. Si on recolle bien le feuilletage est aussi globalement compatible et donc $\ff$ peut \^etre rendu totalement g\'eod\'esique.
$\Box$
\hip 
Pour conclure ce paragraphe notons deux questions qui viennent naturellement:
 \begin{itemize}
\item existe-t-il des vari\'et\'es poss\'edant un flot riemannien mais ne poss\'edant pas de  feuilletage totalement g\'eod\'esique de codimension $1$ ?
\item existe-t-il des vari\'et\'es ne  poss\'edant pas de flot riemannien mais  poss\'edant  des  feuilletages totalement g\'eod\'esiques de codimension $1$ (dont les feuilles de type lumi\`ere sont isol\'ees) ?
\end{itemize}
\subsection{Un resultat sans l'hypoth\`ese <<flot riemannien>>.}
On a donc montr\'e que \emph{toute vari\'et\'e ferm\'ee de dimension $3$ poss\'edant un flot transversalement orientable riemannien poss\`ede aussi un feuilletage totalement g\'eod\'esique}. La proposition suivante est une sorte de r\'eciproque sous certaines hypoth\`eses, essentiellement l'abscence de feuilles de type temps.
\begin{prop}
Soit $M$ une vari\'et\'e ferm\'ee de dimension $3$ munie d'un feuilletage $\ff$  transversalement orientable de codimension $1$ et totalement g\'eod\'esique dont les feuilles de type lumi\`ere sont compactes et en nombre fini non nul et n'ayant pas de feuilles de type temps. Alors $M$ est un fibr\'e de Seifert ou un fibr\'e en tore sur le cercle.
\end{prop}
{\bf Id\'ee de preuve.}  Soit $U$ une  composante connexe  de $M$ priv\'e des feuilles de type lumi\`ere, $U$ poss\'ede un flot riemannien $\phi_U$. Les composantes connexes du bord de $U$ sont forc\'ement des tores. Il y a alors deux possibilit\'es.
Soit toutes les feuilles de $\phi$ vont s'accumuler sur le bord et alors $U\simeq T^2\times]0,1[$, soit l'adh\'erence de toutes les trajectoires de $\phi$ sont des cercles ou des tores et alors $U\simeq D^2\times S^1$ ou  $T^2\times]0,1[$ s'il existe des trajectoires non ferm\'ees ou un fibr\'e de Seifert \`a bord si toutes les trajectoires sont ferm\'ees. Les conditions de recollement font de $M$ soit un fibr\'e de Seifert, soit un fibr\'e en tores.$\Box$
\hip 
Pour pouvoir consid\'erer les feuilles de type temps il manque une connaissance des adh\'erences des feuilles d'un flot lorentzien. En fait il s'agirait de  savoir quelles  sont les vari\'et\'es \`a bord poss\'edant un flot riemannien ou lorentzien tangent au bord et riemannien sur le bord. 
\subsection{Une remarque sur la compl\'etude g\'eod\'esique.}
Nous nous posons la question de la compl\'etude g\'eod\'esique des m\'etriques rendant ces feuilletages totalement g\'eod\'esiques. Nous aimerions reprendre les techniques utilis\'ees par l'auteur dans \cite{SMF}. Cependant celles-ci n\'ecessitent la pr\'esence d'un feuilletage g\'eod\'esique de dimension $1$. Or ici rien ne dit qu'un tel feuilletage soit pr\'esent. Par contre restreint aux feuilles compactes de type lumi\`ere de $\ff$ le flot $\phi$ est  g\'eod\'esique de type lumi\`ere (nul-pr\'eg\'eod\'esique si on reprend la terminologie de \cite {SMF}).
Nous allons ici n'\'etudier que  ces g\'eod\'esiques. Le r\'esultat obtenu ne concerne en fait que les g\'eod\'esiques ferm\'ees d'une feuille de type lumi\`ere au voisinage de laquelle les feuilles de $\ff$ ne changent pas de type.
\begin{prop}
Soit $\ff$ un feuilletage totalement g\'eod\'esique de codimension $1$ transversalement orientable sur une vari\'et\'e lorentzienne $(M,g)$. Soit $F$ une feuille de $\ff$ de type lumi\`ere dont les feuilles voisines sont toutes de type espace ou toutes de type temps.
On suppose de plus que $F$ contient une  g\'eod\'esique $c$ de type lumi\`ere ferm\'ee\footnote{la g\'eod\'esique est dite ferm\'ee simplement si sa courbe se referme, son image dans le fibr\'e tangent n'est pas forc\'ement ferm\'ee.}. Cette g\'eod\'esique est compl\`ete si et seulement elle ne porte pas d'holonomie lin\'eaire de $\ff$. 
\end{prop}
\pre 
On commence par consid\'erer, au voisinage de $F$ un champ de vecteurs $Y$  partout non nul, tangent \`a $\ff$ et orthogonal \`a $TF$. Si les feuilles voisines de $F$ sont de type temps, on choisit $Y$ de type temps hors de $F$ (on utilise par exemple une m\'etrique riemannienne auxiliaire).
Ensuite on construit un atlas adapt\'e \`a $\ff$ et \`a $Y$ \cad un syst\`eme de coordonn\'ees $(x_1,x_2, \dots , x_n)$ tel que $T\ff=\mathrm{Vect}(\partial x_1,\dots, \partial x_{n-1})$ et que $\partial x_1=Y$.
En tout point de $F$ une g\'eod\'esique port\'ee par $Y$ s'ecrit $c(t)=(x_1(t),0,\dots, 0)$ et l'\'equation 
d'Euler-Lagrange nous dit que~:
$$\ddot x_1(t)+\Gamma^1_{1,1}(\dot x_1(t))^2=0 ,$$
o\`u $\Gamma ^1_{1,1}$ est un symbole de Christoffel. On note $g_{i\, j}$ les coefficients de la matrice de $g$ dans ce syt\`eme de coordonn\'ees et $g^{i\,j}$ ceux de la matrice inverse.
On a $$\Gamma ^1_{1,1}=1/2 g^{1\,n}\left (2\partial _1 g_{1\,n} + \partial _n g_{1\,1}\right).$$
C'est ici qu'intervient l'hypoth\`ese concernant le type des feuilles voisines de $F$.  On voit que $g_{1\,1}$ atteint son minimum (ou maximum) sur $F$ et donc $\partial _n g_{1\,1}=0$. Sur la feuille $F$ on a de plus $g^{1\,n}=1/g_{1\,n}$.
L'\'equation d'Euler-Lagrange devient
$$g_{1\,n}\ddot x_1(t)+ \partial_1 g_{1\,n} (\dot x_1(t))^2=0$$
\cad $g_{1\,n} \dot x_1(t)=C$. La constante $C$ d\'epend uniquement de la vitesse initiale de la g\'eod\'esique.

Une g\'eod\'esique ferm\'ee est compl\`ete si et seulement si lorsqu'on l'a parcourue une fois on revient avec la vitesse initiale (elle est alors dite p\'eriodique). Si $c$ est une courbe ferm\'ee on voit que $c$ est compl\`ete si et seulement si apr\`es un tour on retrouve la m\^eme constante $C$. Le fait d'avoir choisi les coordonn\'ees de telle sorte que $\partial x_1=Y$ permet de "suivre" les constantes. Si l'on passe d'un ouvert de l'atlas $U_j$ \`a un ouvert $U_i$ par un changement de carte $\psi^{i,j}$, on relie les constantes $C_i$ et $C_j$ par $C_i=C_j \, \partial_n \psi^{i,j}_n$, o\`u $\partial_n \psi^{i,j}_n$ d\'esigne la d\'eriv\'ee par rapport \`a la n$^\mathrm{ieme}$ coordonn\'ee de la n$^\mathrm{ieme}$ composante du changement de carte, \cad de l'holonomie lin\'eaire du feuilletage $\ff$. Si apr\`es un nombre fini de changements de cartes on est revenu au point initial, on voit que la vitesse est la m\^eme si et seulement si le produit des d\'eriv\'ees est \'egal \`a $1$, \cad si le lacet $c$ ne porte pas d'holonomie lin\'eaire.
$\Box$
\hip
On notera que la plupart des exemples construits plus haut ont des feuilles de type lumi\`ere sans holonomie lin\'eaire. Il semble difficile de conjecturer quelquechose sur la compl\'etude de ces m\'etriques en g\'en\'eral en l'\'etat des connaissances. 

\end{document}